\documentclass[12pt]{amsart}
\usepackage{a4wide}
\usepackage{amssymb}
\usepackage{amsthm}
\usepackage{amsmath}
\usepackage{amscd}
\usepackage[mathscr]{eucal}
\usepackage[all]{xy}

\usepackage{verbatim}

\numberwithin{equation}{section}

\theoremstyle{plain}
\newtheorem{theorem}{Theorem}[section]

\newtheorem{lemma}[theorem]{Lemma}
\newtheorem{proposition}[theorem]{Proposition}

\theoremstyle{definition}
\newtheorem{definition}[theorem]{Definition}
\newtheorem{remark}[theorem]{Remark}
\newtheorem{example}[theorem]{Example}

\theoremstyle{remark}

\newcommand{\R}{\mathbb{R}}
\newcommand{\Q}{\mathbb{Q}}
\newcommand{\Z}{\mathbb{Z}}

\newcommand{\C}{\mathbb{C}}

\renewcommand{\H}{\mathbb{H}}
\newcommand{\F}{\mathbb{F}}

\newcommand{\zxz}[4]{\begin{pmatrix} #1 & #2 \\ #3 & #4 \end{pmatrix}}
\newcommand{\abcd}{\zxz{a}{b}{c}{d}}
\newcommand{\leg}[2]{\left( \frac{#1}{#2} \right)}
\newcommand{\kzxz}[4]{\left(\begin{smallmatrix} #1 & #2 \\ #3 & #4\end{smallmatrix}\right) }
\newcommand{\kabcd}{\kzxz{a}{b}{c}{d}}

\newcommand{\calG}{\mathcal{G}}
\newcommand{\calH}{\mathcal{H}}

\newcommand{\calO}{\mathcal{O}}

\newcommand{\calQ}{\mathcal{Q}}

\newcommand{\fraka}{\mathfrak a}
\newcommand{\frakb}{\mathfrak b}

\newcommand{\frake}{\mathfrak e}

\newcommand{\bs}{\backslash}

\newcommand{\Span}{\operatorname{span}}

\newcommand{\Sl}{\operatorname{SL}}
\newcommand{\Gl}{\operatorname{GL}}

\newcommand{\Symp}{\operatorname{Sp}}
\newcommand{\Mat}{\operatorname{M}}

\newcommand{\Orth}{\operatorname{O}}
\newcommand{\Uni}{\operatorname{U}}

\newcommand{\sig}{\operatorname{sig}}

\newcommand{\End}{\operatorname{End}}

\newcommand{\ord}{\operatorname{ord}}

\begin{document}

\title[The Weil representation and Hecke operators]{The Weil representation and Hecke operators for vector valued modular forms}

\author[Jan H.~Bruinier and Oliver Stein]{Jan Hendrik Bruinier and Oliver Stein}

\address{Mathematisches Institut, Universit\"at zu K\"oln, Weyertal 86--90, D-50931 K\"oln, Germany}
\email{bruinier@math.uni-koeln.de } \email{ostein@math.uni-koeln.de
}


\date{\today}

\begin{abstract}
We define Hecke operators on vector valued modular forms
transforming with the Weil representation associated to a
discriminant form. We describe the properties of the corresponding
algebra of Hecke operators and study the action on modular forms.
\end{abstract}

\maketitle

\section{Introduction}

Hecke operators are a fundamental tool in the study of modular
forms. They can be used to obtain information on the arithmetic
nature of the Fourier coefficients. They are vital for the
definition of $L$-functions associated to modular forms and for
understanding their properties. The theory of Hecke operators is
well developed for scalar valued modular forms \cite{Sh}.

In many recent works, vector valued modular forms associated to the
Weil representation play an important role, see e.g.  \cite{Bo1},
\cite{Bo2}, \cite{Br}, \cite{McG}, \cite{Sch}. For instance,
Borcherds uses them to provide a elegant description of the Fourier
expansion of various theta liftings.  The purpose of the present
paper is to develop the foundations of a Hecke theory for such
vector valued modular forms.
The results can be used to
associate an $L$-function to a vector valued modular form,
essentially the standard $L$-function.

We now describe the content of this paper in more detail.  Let $L$ be
a non-degenerate even lattice of type $(b^+,b^-)$ and level $N$.  The
modulo $1$ reduction of the quadratic form on the dual lattice $L'$
defines a $\Q/\Z$-valued quadratic form on the discriminant group
$A=L'/L$. To simplify the exposition, we assume throughout the
introduction that the signature $\sig(L)=b^+-b^-$ of $L$ is {\em
  even}. In the body of the paper, both, odd and even signature is
treated.

The Weil representation associated to the discriminant form $A$ is a
unitary representation of $\Gamma=\Sl_2(\Z)$ on the group ring
$\C[A]$,
\[
\rho_A: \Gamma\longrightarrow \Uni(\C([A])),
\]
defined by \eqref{eq:weilt}, \eqref{eq:weils}.
It factors through the finite quotient
\[
S(N):=\Sl_2(\Z/N\Z)\cong\Gamma/\Gamma(N).
\]

Let $k\in \Z$.
A holomorphic function $f:\H\to \C[A]$ is called a modular form of
weight $k$ and type $\rho_A$ for the group $\Gamma$, if
\[
f(M\tau)= (c\tau+d)^{k} \rho_A(M) f(\tau)
\]
for all $M=\kabcd\in \Gamma$, and $f$ is holomorphic at the cusp
$\infty$. We denote the vector space of such holomorphic modular
forms by $M_{k,A}$, and write $S_{k,A}$ for the subspace of cusp
forms.

In order to define Hecke operators on vector valued modular forms of
type $\rho_A$, we need to extend the representation $\rho_A$ to a
sufficiently large subgroup of $\Gl_2^+(\Q)$.
%
A natural starting point is to try to extend $\rho_A$, viewed as a
representation of $S(N)$, to a
representation of  $G(N):=\Gl_2(\Z/N\Z)$.  However, it
was observed by E.~Freitag that such an extension does {\em not\/}
exist in general, see Example \ref{exfreitag}.

Here we consider the subgroup
\[
\{M\in G(N);\;\det(M)\equiv \square\pmod{N}\}
\]
of matrices whose determinant is a square modulo $N$.
It has the extension
\[
Q(N)=\{(M,r)\in G(N) \times U(N);\quad \det(M)\equiv
r^2\pmod{N}\},
\]
where $U(N)$ denotes the unit group of $\Z/N\Z$.  The group $Q(N)$
is isomorphic to $S(N)\times U(N)$. Consequently, we may extend the
Weil representation to $Q(N)$ by taking the tensor product of
$\rho_A$ on $S(N)$ and a suitable character on $U(N)$, see
Proposition \ref{thm:main1}.

If $M$ is an element of $G(N)$ whose determinant is  a square modulo
$N$, and $r,r'\in U(N)$ with $\det(M)\equiv r^2\equiv
r'{}^2\pmod{N}$, then $(M,r)$ and $(M,r')$ both belong to $Q(N)$. We
prove that the action of $\rho_A(M,r)$ and $\rho_A(M,r')$ on $\C[A]$
differ only by the action of an element of the orthogonal group
$\Orth(A)$, see Proposition \ref{prop:choices}.

This extension of the Weil representation can be used to define a
Hecke operator $T(M,r)$ on $M_{k,A}$ for every pair $(M,r)$, where
$M\in \Mat_2(\Z)$ and $r\in U(N)$ with $\det(M)\equiv r^2\pmod{N}$.
We compute the action of these operators on the Fourier expansion of a
modular form (see Section \ref{subsect:4.1}).  They
  generalize the classical Hecke operators on scalar valued
modular forms and Jacobi forms
(see e.g. Remark \ref{compareBB} and Remark \ref{compareJacobi}).

In particular, for every positive integer $m$ coprime to $N$ we
obtain a Hecke operator
\[
T(m^2)^*:=T\left(\kzxz{m^2}{0}{0}{1},m\right)
\]
on $M_{k,A}$. These operators generate a commutative subalgebra of
$\End(M_{k,A})$, which is actually already generated by the
$T(p^2)^*$ for $p$ prime and coprime to $N$. The operators
$T(m^2)^*$ take cusp forms to cusp forms and are self-adjoint with
respect to the Petersson scalar product (see Theorem
\ref{thm:structure1}). In particular, $S_{k,A}$ has a basis of
simultaneous eigenforms of all $T(m^2)^*$ with $(m,N)=1$.

In Section \ref{sect:5} we extend the definition of the Hecke
operators $T(m^2)^*$ to all positive integers $m$, not necessarily
coprime to $N$. This is done by defining the right-action on $\C[A]$
of a matrix $\alpha=\kzxz{m^2}{0}{0}{1}$  by the same formula as in
the case where $m$ is coprime to $N$. Notice, that the corresponding
linear map
\[
\C[A]\longrightarrow \C[A],\quad \frake_\lambda\mapsto
\frake_\lambda\mid_A\alpha=\frake_{m\lambda}
\]
is neither surjective nor injective in general.  However, it still
can be used to obtain an ``action'' of the double coset $\Gamma
\alpha \Gamma $, see Proposition \ref{welldef} and Lemma
\ref{lem:cosetact}. This suffices to define a corresponding Hecke
operator $T(m^2)^*$ on $M_{k,A}$, which is consistent with our
earlier definition when $m$ is coprime to $N$.

For any positive integer $m$, the Hecke operator $T(m^2)^*$ is self
adjoint with respect to the Petersson scalar product. Moreover, if
$m$ and $n$ are coprime, then
\[
T(m^2)^*T(n^2)^* = T(m^2 n^2)^*,
\]
see Theorem \ref{thm:structure2}.  Observe that for a prime $p$
dividing $N$ the local Hecke algebra, that is, the subalgebra of
$\End(M_{k,A})$ generated by the $T(p^{2\nu})^*$, is considerably more
complicated than in the case where $p$ is coprime to $N$. For
instance, it is commutative if $p$ is coprime to $N$, but in general
non-commutative if $p$ divides $N$.

Let $S$ be a finite set of primes and let $N_S$ be the product of the
primes in $S$.
Let $f\in S_{k,A}$ be a common eigenform of all $T(m^2)^*$ with $(m,N_S)=1$, so
\[
f\mid_{k,A} T(m^2)^*= \lambda_m(f) f.
\]
We can use the above results to define an $L$-function
associated to $f$ by putting
\[
L^S(s,f)= \sum_{\substack{m\geq 1\\ (m,N_S)=1}} \lambda_m(f) m^{-s}.
\]
It is easily seen that $L^S(s,f)$ converges for $\Re(s)$
sufficiently large.  By Theorem \ref{thm:structure2}, this
$L$-function has an Euler product expansion. According to
\cite{Boe}, it should be viewed as the standard $L$-function of $f$.
It would be interesting to study the analytic properties of
$L^S(s,f)$ in more detail. This could possibly be done by using a
variant of the doubling method (see \cite{Boe}, \cite{Ga},
\cite{PSR}) involving a Siegel Eisenstein series of genus $2$
associated to the Weil representation of $\Symp(2,\Z)$ on $\C[A^2]$.

Moreover, it would be very interesting to develop a theory of new
forms for the space $M_{k,A}$. One could try to associate an
irreducible automorphic representation to a vector valued new form
and study the properties of the resulting map.

If the signature of $L$ is odd, one can carry over the above
results. However, one has to work with the metaplectic cover of
$\Gamma$ and with similar $\{\pm 1\}$-extensions of $S(N)$, $G(N)$,
and $Q(N)$. In this case $M_{k,A}$ vanishes unless $k$ is
half-integral.  Following the argument of Shimura \cite{Sh2}, we
show that the Hecke operator $T(M,r)$ vanishes identically unless
$\det(M)$ is the square of a rational number, see Proposition
\ref{prop:heckevanish}. The computation of the action of the Hecke
operators on modular forms is more involved than in the even signature
(i.e.  integral weight) case, see Theorem \ref{fourierodd}.


We thank E. Freitag for many valuable discussions on this paper. Moreover, we thank J. Funke for several useful comments.

\section{Discriminant forms and the Weil representation}

Here we briefly summarize some facts on lattices, discriminant
forms, and the Weil representation. See also \cite{Bo1}, \cite{Bo2},
\cite{Br}.

Let $L$ be a non-degenerate even lattice of type $(b^+,b^-)$.
We denote the bilinear form on $L$ by $(\cdot,\cdot)$ and the
associated quadratic form by $x\mapsto
\frac{1}{2}x^2=\frac{1}{2}(x,x)$. We let $\sig(L)=b^+-b^-$ be the
signature of $L$.  We write $L'$ for the dual lattice of
$L$, and denote by $N$ the level of $L$, that is, the smallest
positive integer such that $\frac{N}{2}x^2\in \Z$ for all $x\in L'$.
The finite abelian group $L'/L$ is called the discriminant group of
$L$. Its order is equal to the absolute value of the Gram
determinant of $L$.

Recall that a discriminant form is a finite abelian group $A$
together with a $\Q/\Z$-valued non-degenerate quadratic form
$x\mapsto \frac{1}{2}x^2$, for $x\in A$ (see \cite{Ni}).  If $L$ is
a non-degenerate even lattice then $L'/L$ is a discriminant form
where the quadratic form is given by the mod $1$ reduction of the
quadratic form on $L'$. Conversely, every discriminant form can be
obtained in this way.  The quadratic form on $L'/L$ determines the
signature of $L$ modulo $8$ by Milgram's formula (see \cite{MH}
Appendix 4):
\begin{equation}
\label{eq:milgram} \sum_{\lambda\in L'/L}
e(\lambda^2/2)=\sqrt{|L'/L|} e( \sig(L)/8).
\end{equation}
Here and throughout we abbreviate $e(z)=e^{2\pi i z}$ for $z\in \C$.
We define the signature $\sig(A)\in \Z/8\Z$ of a discriminant form
$A$ to be the signature of any even lattice with that discriminant
form.

If $A$ is a discriminant form, then we write $A^n$ for the subgroup
of elements that are $n$-th powers of elements of $A$. Moreover, we
write $A_n$ for the subgroup of elements of $A$ whose order divides
$n$. We have an exact sequence
\begin{align}
\label{eq;sequence} 0\longrightarrow A_n \longrightarrow A
\longrightarrow A^n \longrightarrow 0,
\end{align}
and $A^n$ is the orthogonal complement of $A_n$.

Let $\H=\{\tau\in \C;\;\Im(\tau)>0\}$ by the complex upper half
plane. We write $\widetilde{\Gl}_2^+(\R)$ for the metaplectic
two-fold cover of $\Gl_2^+(\R)$. The elements of this group are
pairs $(M,\phi(\tau))$ where $M=\kabcd\in\Gl_2^+(\R)$ and
$\phi:\H\to \C$ is a holomorphic function with
$\phi(\tau)^2=c\tau+d$.  The multiplication is defined by
\[
(M,\phi(\tau)) (M',\phi'(\tau))=(M M',\phi(M'\tau)\phi'(\tau)).
\]
For $g=(M,\phi)\in \widetilde{\Gl}_2^+(\R)$, we put
$\det(g)=\det(M)$.
Moreover, if $G$ is a subset of $\Gl_2^+(\R)$, we write $\tilde{G}$
for its inverse image under the covering map. Throughout we write
$\Gamma=\Sl_2(\Z)$ for the full modular group. It is well known that
the integral metaplectic group $\tilde\Gamma$ is generated by $T=
\left( \kzxz{1}{1}{0}{1}, 1\right)$, and $S= \left(
\kzxz{0}{-1}{1}{0}, \sqrt{\tau}\right)$. One has the relations
$S^2=(ST)^3=Z$, where $Z=\left( \kzxz{-1}{0}{0}{-1}, i\right)$ is
the standard generator of the center of $\tilde\Gamma$.

We now recall the Weil representation associated with a discriminant
form $A$ (see also \cite{Bo1}, \cite{Bo2}). It is a representation
of $\tilde\Gamma$ on the group algebra $\C[A]$.  We denote the
standard basis elements of $\C[A]$ by $\frake_\lambda$, $\lambda\in
A$, and write $\langle\cdot,\cdot \rangle$ for the standard scalar
product (antilinear in the second entry) such that $\langle
\frake_\lambda,\frake_\mu\rangle =\delta_{\lambda,\mu}$. The Weil
representation $\rho_A$ associated with the discriminant form $A$ is
the unitary representation of $\tilde\Gamma$ on the group algebra
$\C[A]$ defined by
\begin{align}
\label{eq:weilt}
\rho_A(T)(\frake_\lambda)&=e(\lambda^2/2)\frake_\lambda,\\
\label{eq:weils}
\rho_A(S)(\frake_\lambda)&=\frac{e(-\sig(A)/8)}{\sqrt{|A|}} \sum_{\mu\in A} e(-(\lambda,\mu)) \frake_\mu.\\
\intertext{Note that}
\label{eq:weilz} \rho_A(Z)(\frake_\lambda)&=e(-\sig(A)/4)
\frake_{-\lambda}.
\end{align}
The orthogonal group $\Orth(A)$ also acts on $\C[A]$ by
\begin{align}
\label{eq:weilh} \rho_A(h)(\frake_\lambda)&=\frake_{h^{-1}\lambda}
\end{align}
for $h\in\Orth(A)$, and the actions of $\tilde\Gamma$ and $\Orth(A)$
commute.

If the signature of $A$ is even, then \eqref{eq:weilz} implies that
$Z^2$ acts trivially. Hence, the Weil representation factors through
$\Gamma$. Moreover, it is trivial on the principal congruence
subgroup $\Gamma(N)$, where $N$ is the level of $A$, i.e., the level
of any even lattice $L$ with $L'/L=A$ (see e.g.~\cite{Eb}, Chapter
3, Theorem~3.2). Therefore, $\rho_A$ factors through the finite
group
\begin{align}
\label{quotev} S(N):=\Sl_2(\Z/N\Z)\cong\Gamma/\Gamma(N).
\end{align}
If the signature of $A$ is odd, we notice that the level of $A$ must
be divisible by $4$. This follows from the oddity formula (\cite{CS}
p.~383 (30)) which implies that $A$ contains odd $2$-adic Jordan
components. On $\Gamma(4)$ the metaplectic cover has the section
\[
s:\Gamma(4)\longrightarrow \tilde{\Gamma}(4),\quad
\abcd\mapsto\left( \abcd, \leg{c}{d}\sqrt{c\tau+d}\right)
\]
given by the theta multiplier system. Here $\sqrt{\cdot}$ denotes
the principal branch of the holomorphic square root. The same
argument as at the end of the proof of Theorem 5.4 in \cite{Bo2}
implies that $\rho_A$ is trivial on $s(\Gamma(N))$ and factors
through the central extension of $S(N)$ by $\{\pm 1\}$ given by
\begin{align}
\label{quotodd} S_1(N):=\tilde{\Gamma}/s(\Gamma(N)).
\end{align}



We will also need the action of $\rho_A$ on diagonal matrices in
$S(N)$. Following \cite{McG}, for integers $a,d$ coprime to $N$ such
that $ad\equiv 1\pmod {N}$, we put
\begin{align}
R_d:=ST^dS^{-1}T^a ST^d.
\end{align}
It is easily checked that $R_d=(M,\phi)$ where $M\equiv
\kzxz{a}{0}{0}{d}\pmod{N}$.

\begin{lemma}
\label{mcg4.6} (See \cite{McG} Lemma 4.6.) For $a,d$ as above we
have
\begin{align}
\label{eq:2.10}
 \rho_A(R_d)\frake_\lambda = \frac{g_d(A)}{g(A)}
\frake_{d\lambda}.
\end{align}
Here $g_d(A)$ denotes the Gauss sum
\begin{align}
g_d(A)= \sum_{\lambda\in A} e(d\lambda^2/2)
\end{align}
and $g(A)=g_1(A)$. \hfill $\square$
\end{lemma}

Notice that by Milgram's formula we have $g(A)= \sqrt{|A|}
e(\sig(A)/8)$.  Moreover, one easily checks that
$|g_d(A)|=\sqrt{|A|}$. If $r\in \Z$ is coprime to $N$, then we have
$g_{dr^2}(A)=g_d(A)$.  In particular, $g_d(A)=g_a(A)$.  Finally,
Lemma \ref{mcg4.6} and the fact that $\rho_A$ is a representation
imply the relation
\begin{align}
\label{gausscocycle}
\frac{g_{dr}(A)g(A)}{g_d(A)g_r(A)}=\begin{cases}1,&\text{if $\sig(A)$ is even,}\\
\pm 1,&\text{if $\sig(A)$ is odd.}
                   \end{cases}
\end{align}

The following more general formula was given by Borcherds.

\begin{proposition} \label{bo-thm5.4}
(See \cite{Bo2} Theorem 5.4.) Let $g=(\kabcd, \sqrt{c\tau+d})\in
\tilde\Gamma$, and suppose that $b$ and $c$ are divisible by $N$.
Then
\begin{align}
\rho_A(g)\frake_\lambda= \chi_A(g) \frake_{d\lambda}.
\end{align}
Here $\chi_A$ denotes the character of $\tilde\Gamma_0(N)$ defined
in \cite{Bo2} Theorem 5.4.\hfill $\square$
\end{proposition}

\begin{lemma}
\label{uact} Let $U=\left(\kzxz{1}{0}{1}{1},\sqrt{\tau+1}\right)\in
\tilde \Gamma$. The action of $U^m$ is given by
\[
\rho_A(U^m)\frake_\lambda =\frac{1}{|A|}\sum_{\mu,\nu\in A} e\left(
-m\mu^2/2+(\mu,\lambda-\nu)\right) \frake_\nu.
\]
\end{lemma}

\begin{proof}
Since $U^m=ST^{-m}S^{-1}$, this follows immediately from
\eqref{eq:weilt} and \eqref{eq:weils}.
\end{proof}

In many recent works vector valued modular forms associated to the
Weil representation are considered (see e.g. \cite{Bo1}, \cite{Bo2},
\cite{Br}, \cite{McG}, \cite{Sch}).  Let $k\in \frac{1}{2}\Z$, and let
$\Gamma'\subset \tilde\Gamma$ be a subgroup of finite index.
A holomorphic function $f:\H\to \C[A]$ is called a modular form of
weight $k$ and type $\rho_A$ for the group $\Gamma'$, if
\[
f(M\tau)= \phi(\tau)^{2k} \rho_A(M,\phi) f(\tau)
\]
for all $(M,\phi)\in \Gamma'$, and $f$ is holomorphic at the
cusps of $\Gamma'$. We denote the $\C$-vector space of such holomorphic
modular forms by $M_{k,A}(\Gamma')$. Moreover, for the full modular group
we put
$M_{k,A}=M_{k,A}(\tilde\Gamma)$.
Formula \eqref{eq:weilz} implies that
$M_{k,A}=\{0\}$ unless
\begin{align}
\label{eq:parity} 2k\equiv \sig(A)\mod 2.
\end{align}

Recall that for $f,g\in M_{k,A}(\Gamma')$  the
Petersson scalar product is defined by
\begin{align}
(f,g)= \frac{1}{[\tilde\Gamma:\Gamma']}\int_{\Gamma'\bs \H}  \langle f(\tau), g(\tau)\rangle \,y^k\,\frac{dx\,dy}{y^2}.
\end{align}
Here $x$ denotes the real part and $y$ the imaginary part of $\tau\in \H$.
The Petersson scalar product converges when $ f\otimes g$ is a cusp form.

\section{Extending the Weil representation}
\label{sect:4}

In the classical theory of {\em scalar valued} modular forms Hecke
operators play an important role (see e.g.~\cite{Sh}). It is natural
to try to define Hecke operators on vector valued modular forms of
type $\rho_A$ as well.  This requires the extension of the
representation $\rho_A$ to a representation (of a sufficiently large
subgroup) of $\widetilde{\Gl}_2^+(\Q)$.  However, it is not obvious
how this can be done.

A natural starting point is to try to extend $\rho_A$, viewed as a
representation of $S(N)$ (respectively $S_1(N)$), to a
representation of (a double cover of) $\Gl_2(\Z/N\Z)$.  However, it
was observed by E.~Freitag that such an extension does {\em not\/}
exist in general. This follows from the following example.

\begin{example}
\label{exfreitag}
  Let $d\equiv 1 \pmod 4$ be an integer such that $p:=|d|$ is a prime.
  We consider the ring of integers $\calO$ in the quadratic field
  $\Q(\sqrt{d})$ of discriminant $d$. Together with the norm form, it
  is an even lattice of type $(1,1)$ if $d>0$, and of type
  $(2,0)$ if $d<0$. The dual lattice is $\frac{1}{\sqrt{d}}\calO$, the
  inverse of the different, and the corresponding discriminant form
  $A$ can be identified with the finite field $\F_p$ together with the
  quadratic form $x\mapsto -\frac{1}{d}x^2$.  The associated Weil
  representation $\rho_A$ is a $p$-dimensional representation of
  $S(p)=\Sl_2(\F_p)$ on $\C[A]$. The action of the orthogonal group
  $\Orth(A)=\{\pm 1\}$ splits $\C[A]$ into two $S(p)$-invariant
  subspaces
\begin{align*}
\C[A]^+&=\Span\{\frake_\lambda+\frake_{-\lambda};\; \lambda\in A\},\\
\C[A]^-&=\Span\{\frake_\lambda-\frake_{-\lambda};\; \lambda\in A\}.
\end{align*}
They have dimension $\frac{p+1}{2}$, and $\frac{p-1}{2}$,
respectively.  It follows from \cite{NW}, Theorem 4, that the
corresponding representations of $S(p)$ are irreducible.

On the other hand, the character table of $\Gl_2(\F_p)$ is well
known, see e.g. \cite{FH} \S5.2. It has $p-1$ one-dimensional
representations, $p-1$ irreducible $p$-dimensional representations,
$(p-1)(p-2)/2$ irreducible $(p+1)$-dimensional representations, and
$(p^2-p)/2$ irreducible $(p-1)$-dimsional representations.

Now assume that $p\geq 5$ and that $\rho_A$ has an extension
$\tilde\rho_A$ to a representation of $ \Gl_2(\F_p)$. Because of the
irreducibility of $\C[A]^\pm$, such an extension would have to be a
$p$-dimensional irreducible representation of $\Gl_2(\F_p)$. But
these representations remain irreducible under restriction to
$\Sl_2(\F_p)$, see \cite{FH}, p.72 (2). We obtain a contradiction.
\end{example}

\begin{remark}
In \cite{McG}, McGraw continues $\rho_A$ to an action of
$\Gl_2(\Z/N\Z)$. However, this action is {\em not\/} $\C$-linear,
causing serious difficulties when one tries to define Hecke
operators.
\end{remark}

Here we consider a different group extension $Q(N)$
of $S(N)$ and show that $\rho_A$ can be continued to a
representation of $Q(N)$. Together with the considerations of
Section \ref{sect:5} this will suffice for many applications of
Hecke operators; for instance, to define the standard $L$-function
of a modular form of type $\rho_A$.

Let $A$ be a discriminant form as in the previous section, and let
$N$ be the level of $A$.  We denote by $U(N)$ the unit group of
$\Z/N\Z$.
%
We briefly write $G(N)=\Gl_2(\Z/N\Z)$ for the general linear group
modulo $N$. The determinant homomorphism $G(N)\to U(N)$ gives rise
to the exact sequence
\begin{align}
\label{seq1}
 1\longrightarrow S(N)\longrightarrow G(N)\longrightarrow U(N) \longrightarrow 1.
\end{align}
This sequence splits and $G(N)$ can be viewed as a semidirect
product of $S(N)$ and $U(N)$.

\subsection{The case of even signature}

Throughout this subsection we assume that $\sig(A)$ is even.
Let $Q(N)$ be the group
\begin{align}
\label{eq:qn}
Q(N)=\{(M,r)\in G(N) \times U(N);\quad \det(M)\equiv
r^2\pmod{N}\}
\end{align}
with the product defined component-wise. We have an exact sequence
\begin{align}
\label{seq2} 1\longrightarrow S(N)\longrightarrow
Q(N)\longrightarrow U(N) \longrightarrow 1,
\end{align}
where $S(N)\to Q(N)$ is given by $M\mapsto (M,1)$, and $Q(N)\to
U(N)$ is given by $(M,r) \mapsto r$. The latter homomorphism has the
section
\begin{align}
\label{eq:u}
U(N)\longrightarrow Q(N),\quad r \mapsto \left( \zxz{r}{0}{0}{r},
r\right).
\end{align}
For $(M,r)\in Q(N)$ the assignment $(M,r)\mapsto (M
\kzxz{r}{0}{0}{r}^{-1}, r)$ defines an isomorphism $Q(N)\cong
S(N)\times U(N)$.

We consider the action of $S(N)$ on $\C[A]$ by the Weil
representation $\rho_A$. In view of \eqref{gausscocycle}, the
assignment $r\mapsto  \frac{g(A)}{g_r(A)}$ defines a character of
$U(N)$. We define a unitary representation of $U(N)$ on $\C[A]$ by
putting
\begin{align}
\label{eq:weilr} \rho_A(r)\frake_\lambda = \rho_A\!\left(
\zxz{r}{0}{0}{r}, r\right)\frake_\lambda=\frac{g(A)}{g_r(A)}
\frake_{\lambda}
\end{align}
for $r\in U(N)$.

\begin{proposition}
\label{thm:main1} The Weil representation $\rho_A$ of the group
$S(N)$ extends to a unitary representation of $Q(N)$ by
\eqref{eq:weilr}.
\end{proposition}

\begin{proof}
Since $Q(N)=S(N)\times U(N)$ we only have to check that the actions
of $S(N)$ and $U(N)$ commute. This is obvious.
\end{proof}

\begin{remark}
  Clearly we could take any character of $U(N)$ to define the action
  of the Weil representation on $U(N)$. The above choice is compatible
  with the definition of the Weil representation on double cosets in
  Section \ref{sect:5}. Moreover, it is compatible with the usual
  Hecke operators on scalar valued modular forms, see Remark
  \ref{compareBB} and Remark \ref{compareJacobi}.
\end{remark}

The following proposition shows that the first entry of an element
$(M,r)\in Q(N)$ gives the ``essential contribution'' to the Weil
representation. If we fix $M$ then for the different choices of
$r\in U(N)$ the Weil representation $\rho_A(M,r)$ differs by the
action of an element of the orthogonal group $\Orth(A)$.

\begin{proposition}
\label{prop:choices} Let $(M,r_1), (M,r_2)\in Q(N)$. Then $h:A\to
A$, $\lambda\mapsto r_1 r_2^{-1}\lambda$ is an orthogonal
transformation in $\Orth(A)$, and
\[
\rho_A( M,r_1)\frake_\lambda = \rho_A(M,r_2)\rho_A(h) \frake_\lambda
=\rho_A(M,r_2)\frake_{h^{-1}\lambda}.
\]
\end{proposition}

\begin{proof}
We have $(M,r_2)^{-1}(M,r_1) = (1,r_1r_2^{-1})$, and $t=r_1
r_2^{-1}\in U(N)$ has the property that $t^2\equiv 1\pmod{N}$.
According to \eqref{eq:2.10} and \eqref{eq:weilr}, the action of
$(1,t)$ is given by
\[
\rho_A(1,t)\frake_\lambda = \frake_{t^{-1}\lambda}.
\]
Since $(t\lambda,t\mu) =(t^2\lambda,\mu)=(\lambda,\mu)$ for all
$\lambda,\mu\in A$, multiplication by $t^{-1}$ is an orthogonal
transformation.
\end{proof}

\begin{lemma}
\label{lem:3.4}
 We have
\begin{align}
\label{alphamact}
 \rho_A\left( \zxz{r^2}{0}{0}{1}, r\right) \frake_\lambda
&=  \frake_{r^{-1}\lambda},\\
\rho_A\left( \zxz{1}{0}{0}{r^2}, r\right) \frake_\lambda &=
\frake_{r\lambda}.
\end{align}
\end{lemma}

\begin{proof}
In $Q(N)$ we have
\[
\left( \zxz{r^2}{0}{0}{1}, r\right) = \left(R_{r^{-1}},1\right)
\left( \zxz{r}{0}{0}{r}, r\right).
\]
Therefore the fist formula follows from Lemma \ref{mcg4.6} and
\eqref{eq:weilr}. The second formula follows similarly.
\end{proof}

\subsection{The case of odd signature}
\label{sub:3.2}

Throughout this subsection we assume that $\sig(A)$ is odd.  In this
case the argument of the previous section only shows that the Weil
representation extends to a projective representation of $Q(N)$.
More precisely, it extends to a group homomorphism
\[
Q(N)\longrightarrow \Uni(\C[A])/\{\pm 1\}.
\]
This projective representation gives rise to a 2-cocycle of $Q(N)$
with values in $\{\pm 1\}$. The cocycle defines a central extension
$Q_1(N)$ of $Q(N)$ by $\{\pm 1\}$. The group $S_1(N)$, see
\eqref{quotodd}, can be identified with a subgroup of $Q_1(N)$ (see
Section \ref{sub:4.2} for more details).

\begin{proposition}
\label{thm:main2} The assignment \eqref{eq:weilr} defines an
extension of the Weil representation to a projective representation
of $Q(N)$. It lifts to a unitary representation of $Q_1(N)$.
\end{proposition}

\section{Hecke operators on vector valued modular forms}
\label{sect5}

We now use the results of Section \ref{sect:4} to define Hecke
operators on vector valued modular forms of type $\rho_A$.  We
consider the groups
\begin{align}
\calG(N)&=\{ M\in \Gl_2^+(\Q);\quad \text{$\exists n\in \Z$ with
$(n,N)=1$ such that $nM\in \Mat_2(\Z)$}\\
\nonumber
&\phantom{==============} \text{and $(\det(nM),N)=1$}\},\\
\calQ(N)&=\{ (M,r)\in \calG(N)\times U(N);\quad \text{$\det(M)\equiv
r^2 \pmod{N}$}\}.
\end{align}
We view $\Gamma$ as a subgroup of $\calQ(N)$ by the embedding
$M\mapsto (M,1)$.  We are interested in the action of the Hecke
algebra of the pair $(\calQ(N),\Gamma)$ on modular forms of type
$\rho_A$ and weight $k$.  We have to distinguish the cases whether
$\sig(A)$ is even or odd.

\subsection{The case of even signature}
\label{subsect:4.1}

The composition of the reduction map $\calQ(N)\to Q(N)$ with the
Weil representation $\rho_A:Q(N)\to \Uni(\C[A])$ induces a unitary
representation of $\calQ(N)$ on $\C[A]$, which we will denote by
$\rho_A$ as well.  This left action induces a corresponding right
action by
\begin{align}
\label{right1} \fraka \mid_A (M,r) = \rho_A(M,r)^{-1} \fraka.
\end{align}

In view of \eqref{eq:parity} we only have to consider modular forms
of integral weight $k\in \Z$. Let $f$ be a complex valued function
on $\H$. For $M\in \Gl_2^+(\R)$ the Petersson slash operator is
defined by
\begin{align}
\label{right2} (f\mid_{k} M )(\tau) = \det(M)^{k/2} j(M,\tau)^{-k}
f(M\tau).
\end{align}
This defines a right action of $\calG(N)$ on functions $\H\to \C$.
The center acts by multiplication with $\pm 1$.

If $f:\H\to \C[A]$ is a function we write $f=\sum_{\lambda\in A}
f_\lambda \otimes \frake_\lambda$ for its decomposition in
components with respect to the standard basis of $\C[A]$.  The
tensor product of the above two actions yields a right action of
$\calQ(N)$ on such functions, denoted by
\begin{align}
\label{right3} f\mid_{k,A} (M,r) = \sum_{\lambda\in A} \big(
f_\lambda\mid_k M\big)  \otimes \big(\frake_\lambda\mid_A
(M,r)\big).
\end{align}
Notice that a holomorphic function $f:\H\to \C[A]$ belongs to
$M_{k,A}$, if and only if
\[
f\mid_{k,A} M= f
\]
for all $M\in \Gamma$, and $f$ is holomorphic at the cusp $\infty$.

We consider the Hecke algebra of the pair $(\calQ(N),\Gamma)$ in the
sense of Shimura \cite{Sh}.  If $(M,r)\in \calQ(N)$, the
corresponding double coset decomposes into a finite union of left
cosets
\[
\Gamma \cdot (M,r)\cdot \Gamma = \bigcup_{\gamma\in \Gamma \bs
\Gamma M\Gamma} \Gamma \cdot (\gamma,r).
\]

\begin{definition}
\label{def:defhecke}
For $(M,r)\in \calQ(N)$ we define the corresponding Hecke operator
$T(M,r)$ on $M_{k,A}$ by
\begin{align}
\label{def:hecke} f\mid_{k,A} T(M,r) =
\det(M)^{k/2-1}\sum_{\gamma\in \Gamma \bs \Gamma M \Gamma}
f\mid_{k,A} (\gamma,r),\qquad f\in M_{k,A}.
\end{align}
\end{definition}

The usual argument now shows that $f\mid_{k,A} T(M,r)\in M_{k,A}$.
Hence $T(M,r)$ defines an endomorphism of $M_{k,A}$. A modular form
$f\in M_{k,A}$ has a Fourier expansion of the form
\begin{align}
\label{eq:fourier} f(\tau) = \sum_{\lambda\in A}\sum_{\substack{n\in
\Z+\lambda^2/2}} c(\lambda,n)e(n\tau)\otimes \frake_\lambda.
\end{align}

\begin{theorem} \label{fouriereven1}
Let $p$ be a prime which is a square modulo $N$, and assume that
$r\in U(N)$ with $p\equiv r^2 \pmod{N}$. Let $f\in  M_{k,A}$ and
denote the Fourier expansion as in \eqref{eq:fourier}. Then
\begin{align*}
f|_{k,A} \,T\!\left( \kzxz{p}{0}{0}{1},r\right) = \sum_{\lambda\in
A}\sum_{n\in \Z+\lambda^2/2} b(\lambda,n)e(n\tau)\otimes
\frake_\lambda,
\end{align*}
where
\begin{align*}
b(\lambda,n)=c(r\lambda,pn)+p^{k-1}c(\lambda/r,n/p).
\end{align*}
Here we understand that $c(\lambda/r,n/p)=0$ if $p\nmid n$, i.e., if
$\ord_p(n)=0$.
\end{theorem}

\begin{proof}
Using Lemma \ref{lem:3.4}, the formula follows in the same way as in
the scalar valued case.
\end{proof}

\begin{proposition}\label{fouriereven2}
Let $p$ be a prime coprime to $N$. Let $f\in  M_{k,A}$ and denote
the Fourier expansion as in \eqref{eq:fourier}. Then
\begin{align*}
f|_{k,A} \,T\!\left( \kzxz{p^2}{0}{0}{1},p\right)  =
\sum_{\lambda\in A}\sum_{n\in \Z+\lambda^2/2}
b(\lambda,n)e(n\tau)\otimes \frake_\lambda,
\end{align*}
where
\begin{align*}
b(\lambda,n)=c(p\lambda,p^2n) +\frac{g(A)}{g_p(A)} p^{k-2}
(\delta_p(n)-1) c(\lambda,n) +p^{2k-2}c(\lambda/p,n/p^2),
\end{align*}
and
\begin{align*}
\delta_p(n)=\begin{cases} p,&\text{if $p\mid n$,}\\
0,&\text{if $p\nmid n$.}
\end{cases}
\end{align*}
Moreover, we understand that $c(\lambda/p,n/p^2)=0$ if $p^2\nmid
 n$.
\end{proposition}

\begin{proof}
We omit the proof, which is similar (but easier) to the proof of
Theorem \ref{fourierodd}.
\end{proof}

\begin{remark}
  \label{compareBB} If $|A|=\ell$ is a prime, then $M_{k,A}$ can be
  identified with the plus or minus subspace of the space
  $M_k(\ell,\chi_\ell)$ of scalar valued modular forms for
  $\Gamma_0(\ell)$ with nebentypus, see \cite{BB}. Under this
  identification, the above Hecke operators correspond to the usual
  Hecke operators on $M_k(\ell,\chi_\ell)$. This follows by comparing
  the actions on Fourier expansions. For Hecke operators on
  $M_k(\ell,\chi_\ell)$ see e.g.~\cite{Mi}, Lemma 4.5.14.
\end{remark}

\subsection{The case of odd signature}
\label{sub:4.2}

We now assume that $\sig(A)$ is odd.  In view of \eqref{eq:parity}
we only have to consider modular forms of half integral weight $k\in
\Z+1/2$.  In this case, the problem arises that both, \eqref{right1}
and \eqref{right2}, only define projective actions of $\calQ(N)$. To
obtain honest actions, one has to consider appropriate central
extensions.

We begin by considering the action on $\C[A]$.  The composition of
the natural reduction $\calQ(N)\to Q(N)$ with the projective Weil
representation $\rho_A:Q(N)\to \Uni(\C[A])/\{\pm 1\}$ induces a
projective representation
\[
\rho_A:\calQ(N)\longrightarrow\Uni(\C[A])/\{\pm 1\},\quad g\mapsto
\rho_A (g).
\]
If we choose for every $g\in \calQ(N)$ a $\tilde \rho_A (g)\in
\Uni(\C[A])$ such that $\tilde \rho_A (g)\mapsto \rho_A(g)$ under
the projection to $\Uni(\C[A])/\{\pm 1\}$, we obtain a 2-cocycle $c$
with values in $\{\pm 1\}$ defined by
\[
\tilde \rho_A (g_1 g_2) = c(g_1,g_2)\tilde \rho_A (g_1) \tilde
\rho_A (g_2)
\]
for $g_1,g_2\in \calQ(N)$. This cocycle gives rise to a central
extension
\begin{align}
\label{eq:q2n} \calQ_1(N)= \calQ(N)\times \{\pm 1\},
\end{align}
where the multiplication is defined by
\[
(g_1,t_1)(g_2,t_2) = (g_1 g_2, t_1 t_2 c(g_1,g_2)^{-1}).
\]
Here the cocycle condition
\[
c(g_1 g_2, g_3)c(g_1,g_2)=c(g_1,g_2 g_3)c(g_2,g_3)
\]
is equivalent to the associativity law for the above multiplication.
We obtain a unitary representation
\begin{align}
\label{eq:weilagain} \calQ_1(N)\longrightarrow \Uni(\C[A])
\end{align}
by putting $\rho_A(g,t)=t\tilde\rho_A(g)$. This left action induces
a corresponding right action
\begin{align}
\label{rightodd1} \fraka \mid_A (g,t) = \rho_A(g,t)^{-1} \fraka,
\end{align}
for $(g,t)\in \calQ_1(N)$ and $\fraka\in \C[A]$.

Without loss of generality, for
$(M,1)\in\Gamma\times\{1\}\subset\calQ(N)$ we choose
\begin{align}
\tilde \rho_A (M,1)=\rho_A(M,\sqrt{j(M,\tau)}).
\end{align}
Then we have an injective homomorphism
\begin{align}
\label{hom1} \tilde{\Gamma}\longrightarrow \calQ_1(N),\quad
(M,\pm\sqrt{j(M,\tau)})\mapsto (M,1,\pm 1).
\end{align}
Moreover, for a positive integer $m$ coprime to $N$,  we put
\begin{align}
\label{ximact} \tilde \rho_A \left( \zxz{m^2}{0}{0}{1}, m
\right)\frake_\lambda=\frake_{m^{-1}\lambda}.
\end{align}

To define an action on functions  we consider the
metaplectic group $\widetilde{\Gl}_2^+(\R)$. It acts on functions
$f:\H\to\C$ by
\begin{align}
\label{rightodd2} (f\mid_{k} (M,\phi) )(\tau) = \det(M)^{k/2}
\phi(\tau)^{-2k} f(M\tau)
\end{align}
for $(M,\phi)\in \widetilde{\Gl}_2^+(\R)$. In particular, we have an
action of $\widetilde{\calH}(N)\subset \widetilde{\Gl}_2^+(\R)$,
where
\begin{align}
\label{def:hn} \calH(N)&=\{ M\in \calG(N);\quad \text{$\det(M)$ is a
square mod $N$}\}\subset \Gl_2^+(\R)
\end{align}
is the image of the projection of $\calQ(N)$ to the first component.
Notice that the cocycle of $\calQ(N)$ given by the choice of
$\pm\sqrt{j(M,\tau)}$ for $(M,r)\in \calQ(N)$ and the cocycle $c$
are not isomorphic (this follows for instance from Lemma
\ref{lem:tmn} below). However, their restrictions to $\tilde \Gamma$
are isomorphic.

To define an action on $\C[A]$-valued functions in weight $k$ we
have to consider a combination of the above extensions. We let
$\calQ_{2}(N)$ be the group of tuples $(M,\phi,r,t)$, where
$g=(M,\phi)\in \widetilde{\calH}(N)$, and $r\in U(N)$ with
$\det(M)\equiv r^2\pmod{N}$, and $t\in \{\pm 1\}$. The composition
law is defined by
\begin{align}
\label{eq:q12n} (g_1,r_1,t_1)(g_2,r_2 ,t_2) = (g_1 g_2, \,r_1 r_2,
\,t_1 t_2 c((M_1,r_1),(M_2,r_2))^{-1})
\end{align}
for $(g_i,r_i,t_i)\in \calQ_{2}(N)$ and $g_i=(M_i,\phi_i)$. We
denote by $P:\calQ_{2}\to \widetilde{\calH}(N)$ the natural
projection. It has the kernel
\[
\{(1,1,r,t)\in \calQ_{2};\;r^2\equiv 1 \;(N)\}.
\]
Over $\tilde\Gamma$ the projection $P$ has the section
\begin{align}
L:\tilde\Gamma \longrightarrow \calQ_{2}, \quad
(M,\pm\sqrt{j(M,\tau)})\mapsto (M,\pm \sqrt{j(M,\tau)},1,\pm 1).
\end{align}
We write
\begin{align}
\label{eq:delta} \Delta=L(\tilde\Gamma).
\end{align}
We define the Weil representation $\rho_A$ on $\calQ_{2}(N)$ by
composing the natural map to $\calQ_1(N)$ with the Weil
representation on that group. For $\gamma\in \tilde \Gamma$ we have
\begin{align}
\label{weil12} \rho_A(\gamma)=\rho_A(L(\gamma)).
\end{align}

The tensor product of the above two actions yields a right action of
$\calQ_{2}(N)$ on functions $f:\H\to\C[A]$. If we write
$f=\sum_{\lambda\in A} f_\lambda \otimes \frake_\lambda$ for the
decomposition in components with respect to the standard basis of
$\C[A]$, the action is given by
\begin{align}
\label{right3odd} f\mid_{k,A} (M,\phi,r,t) = \sum_{\lambda\in A}
\big( f_\lambda\mid_k(M,\phi)  \big)  \otimes
\big(\frake_\lambda\mid_A (M,r,t)\big)
\end{align}
for $(M,\phi,r,t)\in \calQ_{2}(N)$. A holomorphic function $f:\H\to
\C[A]$ belongs to $M_{k,A}$, if and only if $f\mid_{k,A} (M,\phi)=
f$ for all $(M,\phi)\in\tilde\Gamma$, and $f$ is holomorphic at the
cusp $\infty$.

We may now define Hecke operators on modular forms of type $\rho_A$
following Shimura \cite{Sh2}. For $\alpha=(M,\phi)\in
\widetilde{\calH}(N)$ and $\xi=(\alpha,r,t)\in\calQ_{2}(N)$ we
consider the double coset $\Delta\xi \Delta$. If
\[
\Delta\xi \Delta = \bigcup_{i} \Delta \xi_i
\]
is a left coset decomposition, we define the Hecke operator $T(\xi)$
by
\begin{align}
\label{heckeodd} f\mid_{k,A} T(\xi) = \det(\alpha)^{k/2-1}\sum_i
f\mid_{k,A} \xi_i
\end{align}
for $f\in M_{k,A}$. It is easily seen that $T(\xi)$ is independent
of the choice of the coset representatives and defines an
endomorphism of $M_{k,A}$. We recall the following standard lemma.


\begin{lemma}
\label{lem:left} Let the notation be as above. Then
\[
\Delta=\bigcup_{i\in I} (\Delta\cap  \xi^{-1}\Delta\xi)\cdot
\gamma_i
\]
($\gamma_i\in \Delta$) is a disjoint left coset decomposition if and
only if
\[
\Delta\xi \Delta = \bigcup_{i\in I} \Delta \cdot \xi\gamma_i
\]
is a disjoint left coset decomposition. \hfill$\square$
\end{lemma}


To compute the action of Hecke operators we have to compare the
groups $\Delta\cap  \xi^{-1}\Delta\xi$ and $L(\tilde\Gamma\cap
\alpha^{-1}\tilde\Gamma\alpha)$. Let $\alpha\in \widetilde{\calH}(N)$ and
$\xi\in \calQ_2(N)$ with $P(\xi)=\alpha$. For $\gamma\in
\tilde\Gamma\cap \alpha^{-1} \tilde\Gamma \alpha$ we have
\begin{align}
L(\alpha \gamma \alpha^{-1})= \xi L(\gamma) \xi^{-1} \cdot
(1,1,t(\gamma))
\end{align}
with $t(\gamma)\in \{\pm 1\}$. Here $t(\gamma)$ is independent of
the choice of $\xi$ and is determined by the condition
\begin{align}
\label{eq:tcond} \rho_A(\alpha \gamma \alpha^{-1}) =
t(\gamma)\rho_A(\xi) \rho_A(\gamma) \rho_A(\xi)^{-1}.
\end{align}
Hence $t$ defines a group homomorphism
\begin{align}
\label{eq:thom} t: \tilde\Gamma\cap \alpha^{-1} \tilde\Gamma
\alpha\longrightarrow \{\pm 1\}.
\end{align}

\begin{lemma}
\label{lem:shimura1} Let the notation be as above. We have
$L(\ker(t))=\Delta\cap  \xi^{-1}\Delta\xi$. Moreover, if $t$ is
non-trivial, then $f\mid_{k,A} T(\xi)=0$ for all $f\in M_{k,A}$.
\end{lemma}

\begin{proof}
The proof is analogous to Proposition 1.0 in \cite{Sh2}
\end{proof}

\begin{lemma}
\label{lem:shimura2} Let the notation be as above. The homomorphism
$t$ is trivial if and only if $P$ gives a bijective map of
$\Delta\xi \Delta$ onto $\tilde \Gamma \alpha \tilde \Gamma$.
Moreover, when this is the case then $\Delta\xi \Delta =\bigcup_{i}
\Delta \xi_i$ (where $\xi_i\in \Delta\xi \Delta$) is a disjoint
union if and only if $\tilde \Gamma \alpha \tilde \Gamma =
\bigcup_{i} \tilde\Gamma P(\xi_i)$ is a disjoint union.
\end{lemma}

\begin{proof}
The proof is analogous to Proposition 1.1 in \cite{Sh2}
\end{proof}

\begin{lemma}
\label{lem:tmn} Let $m,n$ be positive integers coprime to $N$,
$\alpha =\left(\kzxz{m}{0}{0}{n},\sqrt{n}\right)\in
\widetilde{\calH}(N)$, and $\xi=(\alpha,r,t)\in\calQ_{2}(N)$. Define
$t: \tilde\Gamma\cap \alpha^{-1} \tilde\Gamma
\alpha\longrightarrow \{\pm 1\}$ as in \eqref{eq:thom}. Then
\[
t\left( \abcd, \sqrt{c\tau+d}\right) = \leg{mn}{d}.
\]
Here the quadratic residue symbol is defined as in \cite{Sh2}.
\end{lemma}

\begin{proof}
Define $\Gamma'=\Gamma_0(m)\cap \Gamma^0(n)\subset \Gamma$. We first
notice that
\[
\tilde\Gamma'=\tilde\Gamma\cap \alpha^{-1}\tilde\Gamma \alpha.
\]
Since $t$ is a homomorphism it suffices to prove the assertion for a
set of generators of $\tilde\Gamma'$. It is easily verified that
$\tilde\Gamma'$ is generated by $T^n$, $U^m$, and $\tilde
\Gamma'\cap \tilde\Gamma_0^0(N)$.

For $\gamma=T^n\in \tilde\Gamma'$ we have $\alpha \gamma \alpha^{-1}
= T^m$. We compute $t(\gamma)$ using \eqref{eq:tcond}.  There is a
constant $C_\xi$ of modulus $1$ such that
$\rho_A(\xi)\frake_\lambda=C_\xi \frake_{nr^{-1}\lambda}$. Hence
\begin{align*}
\rho_A(\xi) \rho_A(T^n) \rho_A(\xi)^{-1}\frake_\lambda &=
C_\xi^{-1} \rho_A(\xi) \rho_A(T^n) \frake_{n^{-1}r\lambda}\\
&=C_\xi^{-1} e( n^{-1} r^{2} \lambda^2/2)\rho_A(\xi) \frake_{n^{-1}r\lambda}\\
&=e( n^{-1} r^{2} \lambda^2/2)\frake_{\lambda}\\
&=e( m \lambda^2/2)\frake_{\lambda}.
\end{align*}
Here $n^{-1}$ in the exponentials means the inverse of $n$ in
$U(N)$. On the other hand, we have
$\rho_A(T^m)\frake_\lambda=e(m\lambda^2/2)\frake_\lambda$.  Hence
$t(T^n)=1$ in accordance with the formula.

For $\gamma=U^m\in \tilde\Gamma'$ we have $\alpha \gamma
\alpha^{-1}=U^n$. We find
\begin{align*}
\rho_A(\xi) \rho_A(U^m) \rho_A(\xi)^{-1}\frake_\lambda &= C_\xi^{-1}
\rho_A(\xi) \rho_A(U^m) \frake_{n^{-1}r\lambda}\\
&=\frac{C_\xi^{-1}}{|A|} \rho_A(\xi) \sum_{\mu,\nu\in A}
e\left(
-m\mu^2/2+(\mu,n^{-1}r\lambda-\nu)\right) \frake_\nu\\
&= \frac{1}{|A|} \sum_{\mu,\nu\in A} e\left(
-m\mu^2/2+(n^{-1}r\mu,\lambda-\nu)\right) \frake_\nu\\
&= \frac{1}{|A|} \sum_{\mu,\nu\in A} e\left(
-n\mu^2/2+(\mu,\lambda-\nu)\right) \frake_\nu .
\end{align*}
This is equal to $\rho_A(U^n)$. Hence $t(U^m)=1$ in accordance with
the formula.

To compute $t(\gamma)$ for $\gamma\in \tilde \Gamma'\cap
\tilde\Gamma_0^0(N)$, we use the formula for the Weil representation
of Proposition \ref{bo-thm5.4}.
Using the definition one easily checks that if $\gamma \in \tilde
\Gamma'\cap \tilde\Gamma_0^0(N)$, we have
\[
\chi_A(\gamma)=\chi_A(\alpha \gamma \alpha^{-1}) \leg{mn}{d}.
\]
Therefore
\begin{align*}
\rho_A(\xi) \rho_A(\gamma) \rho_A(\xi)^{-1}\frake_\lambda
&= C_\xi^{-1} \rho_A(\xi) \rho_A(\gamma) \frake_{n^{-1}r\lambda}\\
&= \chi_A(\gamma) \frake_{d\lambda}\\
&= \leg{mn}{d} \rho_A(\alpha \gamma \alpha^{-1})\frake_\lambda.
\end{align*}
So $t(\gamma)=\leg{mn}{d}$ as claimed.
\end{proof}

\begin{proposition}
\label{prop:heckevanish}
Let $\alpha=(M,\phi)\in \widetilde{\calH}(N)$ and
$\xi=(\alpha,r,t)\in\calQ_{2}(N)$. Then the Hecke operator $T(\xi)$
on $M_{k,A}$ vanishes identically unless $\det(M)$ is a square in
$\Q$.
\end{proposition}

\begin{proof}
By multiplying with a positive integer we may assume without loss of
generality that $M$ has entries in $\Z$. According to the elementary
divisor theorem for $\Gamma$ we may further assume that
$M=\kzxz{m}{0}{0}{n}$ with positive integers $m,n$. So we may assume
that
\[
\alpha =\left(\zxz{m}{0}{0}{n},\sqrt{n}\right).
\]
Hence the assertion follows from Lemma \ref{lem:shimura1} and Lemma
\ref{lem:tmn}.
\end{proof}

We now study the relation between Hecke operators and Fourier
coefficients. The following theorem is the analogue of Proposition
\ref{fouriereven2} in the odd signature case. It can be viewed as a
generalization of Theorem~1.7 in \cite{Sh2}.

\begin{theorem}\label{fourierodd}
Let $p$ be a prime coprime to $N$, and put
\begin{align*}
\alpha&=\left( \zxz{p^2}{0}{0}{1}, 1\right)\in \widetilde\calH(N),\\
\xi   &=\left(\alpha,p,1\right)\in \calQ_2(N).
\end{align*}
Let $f\in  M_{k,A}$ and write
\begin{align*}
f(\tau) = \sum_{\lambda\in A}\sum_{\substack{n\in \Z+\lambda^2/2}}
c(\lambda,n)e(n\tau)\otimes\frake_\lambda,\\
f|_{k,A} T(\xi) = \sum_{\lambda\in A}\sum_{n\in \Z+\lambda^2/2}
b(\lambda,n)e(n\tau)\otimes \frake_\lambda.
\end{align*}
Then
\begin{align*}
b(\lambda,n)=c(p\lambda,p^2n)
+\epsilon_p^{\sig(A)+\leg{-1}{|A|}}\leg{p}{|A|2^{\sig(A)}} p^{k-3/2}
\leg{-n}{p} c(\lambda,n) +p^{2k-2}c(\lambda/p,n/p^2).
\end{align*}
Here, for an odd integer $d$ we put
\begin{align*}
\epsilon_d=\begin{cases} 1,&\text{if $d\equiv 1\pmod{4}$,}\\
i,&\text{if $d\equiv -1\pmod{4}$.}
\end{cases}
\end{align*}
Moreover, we understand that $c(\lambda/p,n/p^2)=0$ if $p^2\nmid n$.
\end{theorem}

\begin{proof}

To compute $f|_{k,A} T(\xi)$, we need a set of representatives for
$\Delta\backslash \Delta \xi\Delta$. In view of
Lemma~\ref{lem:shimura2} and Lemma~\ref{lem:tmn}, the map
\[
L_\xi: \tilde\Gamma \alpha \tilde \Gamma \longrightarrow \Delta \xi
\Delta, \quad \delta=\gamma\alpha\gamma'\mapsto
L_\xi(\delta):=L(\gamma)\xi L(\gamma')
\]
is a bijection (where $\gamma,\gamma'\in \tilde\Gamma$).  Here we
have $L_\xi(\delta)=(\delta, p,t)$, and $t=t(\delta)$ is uniquely
determined by the condition
\begin{align}\label{weilcond}
\rho_A(\delta,p,t) = \rho_A(\gamma)\rho_A(\xi)\rho_A(\gamma').
\end{align}

We have the disjoint left coset decomposition
\[
\tilde\Gamma \alpha \tilde \Gamma = \tilde \Gamma \alpha \cup
\bigcup_{h\,(p)^*} \tilde \Gamma \beta_h \cup \bigcup_{b\,(p^2)}
\tilde \Gamma \gamma_b,
\]
where
\begin{align*}
\beta_h&=\left( \zxz{p}{hN}{0}{p},\sqrt{p}\right) =\left(
\zxz{r}{Nh}{Ns}{p},\sqrt{Ns\tau +p}\right) \alpha
\left( \zxz{1}{0}{-Nps}{1},\sqrt{-Nps\tau+1}\right),\\
\gamma_b &=\left( \zxz{1}{b}{0}{p^2},p\right)= \left(
\zxz{d}{Nt}{-N}{p^2},\sqrt{-N\tau+p^2}\right)\alpha
\left( \zxz{1}{-Nt}{N}{p^2d},\sqrt{N\tau+p^2d}\right)T^b,\\
\end{align*}
and $r,s\in \Z$ are chosen such that $pr-N^2hs=1$, and $d,t\in \Z$
are chosen such that $p^2 d+N^2t=1$. Consequently, we obtain the
disjoint left coset decomposition
\begin{align}
\Delta \xi \Delta = \Delta \xi \cup \bigcup_{h\,(p)^*} \Delta
L_\xi(\beta_h) \cup \bigcup_{b\,(p^2)} \Delta L_\xi(\gamma_b).
\end{align}

The action of $L_\xi(\beta_h)$ and $L_\xi(\gamma_b)$ in the Weil
representation can be computed by means of \eqref{weilcond} and the
above decompositions.  Using Proposition \ref{bo-thm5.4} and the
fact that $\rho_A(\xi)\frake_\lambda =\frake_{p^{-1}\lambda}$ by
\eqref{ximact}, we find that
\begin{align*}
\rho_A( L_\xi(\gamma_b))\frake_\lambda&=
\chi_A\left(\zxz{d}{Nt}{-N}{p^2},\sqrt{-N\tau+p^2}\right)\\
&\phantom{=}\times\chi_A\left(
\zxz{1}{-Nt}{N}{p^2d},\sqrt{N\tau+p^2d}\right)e(b \lambda^2/2)
\frake_{p\lambda}.
\end{align*}
Since $d$ is a square modulo $N^2$, it is a square modulo the
square-free part of $|A|$ and modulo $8$. Therefore, a quick
calculation shows that the character values are $1$. Consequently,
\begin{align}
\rho_A( L_\xi(\gamma_b))\frake_\lambda=e(b \lambda^2/2)
\frake_{p\lambda}.
\end{align}
In the same way we obtain
\begin{align*}
\rho_A( L_\xi(\beta_h))\frake_\lambda&=
\rho_A\left( \zxz{r}{Nh}{Ns}{p},\sqrt{Ns\tau +p}\right) \rho_A(\xi)\frake_\lambda\\
&= \rho_A\left( \zxz{r}{Nh}{Ns}{p},\sqrt{Ns\tau +p}\right) \frake_{p^{-1}\lambda}\\
&= \chi_A\left( \zxz{r}{Nh}{Ns}{p},\sqrt{Ns\tau +p}\right) \frake_{\lambda}\\
&=
\epsilon_p^{1-\leg{-1}{|A|}-\sig(A)}\leg{Ns}{p}\leg{p}{|A|2^{\sig(A)}}
\frake_{\lambda}
\end{align*}
Here in the last line we have used the explicit formula for
$\chi_A$. Since $-N^2hs\equiv 1\pmod{p}$, we find
\begin{align}
\rho_A(
L_\xi(\beta_h))\frake_\lambda&=\epsilon_p^{1-\leg{-1}{|A|}-\sig(A)}
\leg{-Nh}{p}\leg{p}{|A|2^{\sig(A)}} \frake_{\lambda}.
\end{align}

Now we can compute the Fourier expansion of $f|_{k,A} T(\xi)$. We
have
\begin{align}\label{formulasum}
 f|_{k,A} T(\xi)= p^{k-2} f|_{k,A}\xi + p^{k-2} \sum_{h\,(p)^*} f |_{k,A} L_\xi(\beta_h) + p^{k-2}
\sum_{b\,(p^2)} f |_{k,A} L_\xi(\gamma_b).
\end{align}
For the first summand we find
\begin{align*}
 p^{k-2} f|_{k,A}\xi &=p^{k-2} \sum_{\lambda\in A} \big( f_\lambda\mid_k \alpha  \big)  \otimes \big(\frake_\lambda
\mid_A \xi \big) \\
& =p^{2k-2}  \sum_{\lambda\in A} f_\lambda(p^2\tau)  \otimes
\frake_{p\lambda}.
\end{align*}
For the second summand in \eqref{formulasum} we get
\begin{align*}
&p^{k-2} \sum_{h\,(p)^*} f |_{k,A} L_\xi(\beta_h) =p^{k-2}
\sum_{h\,(p)^*} \sum_{\lambda\in A} \big( f_\lambda\mid_k \beta_h
\big)  \otimes
\big( \frake_\lambda \mid_A L_\xi(\beta_h) \big)\\
&=\epsilon_p^{\sig(A)+\leg{-1}{|A|}-1}\leg{p}{|A|2^{\sig(A)}}
p^{k-2} \sum_{\lambda\in A} \sum_{h\,(p)^*}\leg{-Nh}{p}
f_\lambda(\tau+Nh/p) \otimes \frake_{\lambda}.
\end{align*}
By means of the formula for the Gauss sum
$\sum_{h(p)^*}\leg{h}{p}e(kh/p)=\leg{k}{p}\epsilon_p\sqrt{p}$ we obtain
\[
\sum_{h\,(p)^*}\leg{-Nh}{p} f_\lambda(\tau+Nh/p)=\sqrt{p}\epsilon_p
\sum_{n\in \Z+\lambda^2/2} \leg{-n}{p}c(\lambda,n)e(n\tau),
\]
and therefore
\begin{align*}
&p^{k-2} \sum_{h\,(p)^*} f |_{k,A} L_\xi(\beta_h)\\
&=\epsilon_p^{\sig(A)+\leg{-1}{|A|}}\leg{p}{|A|2^{\sig(A)}}
p^{k-3/2} \sum_{\lambda\in A} \sum_{n\in \Z+\lambda^2/2}
\leg{-n}{p}c(\lambda,n)e(n\tau) \otimes \frake_{\lambda}.
\end{align*}
Finally, for the third summand in \eqref{formulasum} we get
\begin{align*}
p^{k-2} \sum_{b\,(p^2)} f |_{k,A} L_\xi(\gamma_b) &= p^{k-2}
\sum_{\lambda\in A} \sum_{b\,(p^2)}
\big(f_\lambda\mid_k \gamma_b\big)  \otimes \big( \frake_{\lambda}\mid_A L_\xi(\gamma_b)\big)\\
&=p^{-2} \sum_{\lambda\in A} \sum_{b\,(p^2)}
e(-b(p^{-1}\lambda)^2/2) f_\lambda(\tau/p^2+b/p^2)  \otimes
\frake_{p^{-1}\lambda}\\
&=\sum_{\lambda\in A} \sum_{n\in \Z+\lambda^2/2} c(p\lambda,p^2
n)e(n\tau) \otimes \frake_{\lambda}.
\end{align*}
This concludes the proof of the theorem.
\end{proof}

\begin{remark}
\label{compareJacobi}
Let $m$ be a positive integer.
Let $L$ be the lattice $\Z$ with the quadratic form $x\mapsto -mx^2$. Then $L'=\frac{1}{2m}\Z$, and $M_{k,L'/L}$ is isomorphic to the space $J_{k+1/2,m}$ of Jacobi forms of weight $k+1/2$ and index $m$ (cf.~\cite{EZ}, Theorem 5.1).
Under this isomorphism the Hecke operator  $T(\xi)$ of Theorem \ref{fourierodd} corresponds to the Hecke operator $T_p$ on $J_{k+1/2,m}$ defined in \cite{EZ} \S4 (3). This follows from Theorem  \ref{fourierodd} and \cite{EZ}, Theorem 4.5, by comparing the actions on Fourier expansions.
Notice that we have in this particular case
\[
\epsilon_p^{\sig(A)+\leg{-1}{|A|}}\leg{p}{|A|2^{\sig(A)}}=\leg{m}{p}.
\]

For the lattice $L=\Z$ with the quadratic form $ x\mapsto mx^2$ the
space $M_{k,L'/L}$ is isomorphic to the space $J^{skew}_{k+1/2,m}$ of
skew holomorphic Jacobi forms of weight $k+1/2$ and index $m$ defined
in \cite{Sk}. Again, the Hecke operators of Theorem \ref{fourierodd}
correspond to the usual Hecke operators on skew-holomorphic Jacobi
forms.
\end{remark}

\subsection{A Hecke algebra on vector valued modular forms}
\label{sect:4.3}

Using the double coset actions of the previous section, we may
define for every positive integer $m$ coprime to $N$ a Hecke
operator $T(m^2)^*: M_{k,A}\to M_{k,A},$ by
\begin{gather}
\label{hecketm2}
f\mapsto f\mid_{k,A} T(m^2)^* =
\begin{cases}
f|_{k,A} \,T\!\left( \kzxz{m^2}{0}{0}{1},m\right) ,&\text{if $\sig(A)$ is even,}\\
f|_{k,A} \,T\!\left( \kzxz{m^2}{0}{0}{1},1,m,1\right) ,&\text{if
$\sig(A)$ is odd.}
\end{cases}
\end{gather}
In the case of even signature (that is integral weight), the
operator $T(m^2)^*$ differs from the usual Hecke operator $T(m^2)$
which is given by the sum of double cosets consisting of all
integral matrices of determinant $m^2$. This is the reason for our
notation. In the case of odd signature (that is half-integral
weight), the operator $T(m^2)^*$ is analogous to the Hecke operator
in \cite{Sh2} on scalar valued modular forms.

\begin{theorem}
\label{thm:structure1}
The Hecke operators $T(m^2)^*$ (for $m$ coprime to $N$) generate a
commutative subalgebra of $\End(M_{k,A})$, which is actually already
generated by the $T(p^2)^*$ for $p$ prime and coprime to $N$. The
operators $T(m^2)^*$ take cusp forms to cusp forms and are
self-adjoint with respect to the Petersson scalar product.
\end{theorem}

\begin{proof}
Using the actions \eqref{right3} and \eqref{right3odd}, this follows
in the usual way from the properties of the abstract Hecke algebra
of the pair $(\calQ(N),\Gamma)$, respectively $(\calQ_2(N),\Delta)$.
\end{proof}


\section{The Weil representation on double cosets}
\label{sect:5}

We now want to define Hecke operators $T(m^2)^*$ as in Section
\ref{sect:4.3} for {\em all\/} positive integers $m$, not
necessarily coprime to $N$. If $m$ and $N$ are not coprime, then the
reduction of $\kzxz{m^2}{0}{0}{1}$ does not belong to
$\Gl_2(\Z/N\Z)$. So we cannot use the results of the previous
sections. However, it is still possible to extend the Weil
representation $\rho_A$ to the corresponding double coset in a
compatible way, as we will see.

Let $m$ be a positive integer and $\alpha=\left(
\kzxz{m^2}{0}{0}{1},1 \right)\in \widetilde{\Gl}_2^+(\R)$. We define
a right action on $\C[A]$ by
\begin{align}
\label{alphamactgen}
\frake_\lambda \mid_A \alpha = \frake_{m\lambda}.
\end{align}
To lighten the notation, we will frequenty drop the subscript from the
slash operator.  Comparing with \eqref{alphamact} and \eqref{ximact},
we see that \eqref{alphamactgen} is compatible with our earlier
definition in the case that $(m,N)=1$. Moreover, if
$\delta=\gamma\alpha\gamma'\in \tilde\Gamma \alpha \tilde \Gamma$, we
put
\begin{align}
\frake_\lambda \mid \delta = \frake_{\lambda} \mid \gamma \mid
\alpha \mid \gamma'.
\end{align}
We now show that this right action is well defined,
that is, independent of the decomposition of $\delta$.

\begin{proposition}
\label{welldef} Let $\delta
=\gamma\alpha\gamma'=\gamma_1\alpha\gamma_1'\in \tilde\Gamma \alpha
\tilde \Gamma$ (where $\gamma, \gamma', \gamma_1, \gamma_1'\in
\tilde \Gamma$). Then
\[
\frake_{\lambda} \mid \gamma \mid \alpha \mid \gamma' =
\frake_{\lambda} \mid  \gamma_1 \mid \alpha \mid \gamma_1'.
\]
\end{proposition}

\begin{proof}
First, one easily shows that it suffices to prove the proposition in
the case that $\gamma'=\gamma_1=1$. So we have $\delta
=\gamma\alpha=\alpha\gamma_1'$ and need to show that
\begin{align}
\label{toshow} \frake_{\lambda} \mid \gamma \mid \alpha  =
\frake_{\lambda}  \mid \alpha \mid \gamma_1'.
\end{align}
If we write $\gamma=\left(\kabcd,\pm\sqrt{c\tau+d}\right)$, then
$\gamma_1'=\left(\kzxz{a}{b/m^2}{m^2 c}{d},\pm\sqrt{m^2
c\tau+d}\right)$. In particular, $\gamma\in \tilde \Gamma^0(m^2)$
and $\gamma_1'\in \tilde \Gamma_0(m^2)$. It suffices to prove
\eqref{toshow} for $\gamma$ in a set of generators of
$\tilde\Gamma^0(m^2)$.

It is easily seen that $\tilde\Gamma^0(m^2)$ is generated by
$\tilde\Gamma^0(m^2)\cap \tilde\Gamma_0^0(N)$, $T^{m^2}$, and $U$.
For $\gamma\in\tilde\Gamma^0(m^2)\cap \tilde\Gamma_0^0(N)$ the
identity \eqref{toshow} immediately follows from Proposition
\ref{bo-thm5.4}. For $\gamma=T^{m^2}$ it is easily verified as well.

We now consider \eqref{toshow} for $\gamma=U$. Using Lemma
\ref{uact}, we see that the left hand side of \eqref{toshow} is
equal to
\begin{align*}
\frake_{\lambda} \mid U \mid \alpha &=\frac{1}{|A|}\sum_{\mu,\nu\in
A} e\left( \mu^2/2+(\mu,\lambda-\nu)\right) \frake_{m\nu}.
\end{align*}
Using \eqref{eq;sequence}, in the sum over $\nu$ we write
$\nu=\nu'/m+\nu''$ where $\nu'\in A^m$ and $\nu''\in A_m$. We obtain
\begin{align*}
\frake_{\lambda} \mid U \mid \alpha &=\frac{1}{|A|} \sum_{\mu\in A}
\sum_{\substack{\nu'\in A^m\\ \nu''\in A_m}} e\left(
\mu^2/2+(\mu,\lambda-\nu'/m-\nu'')\right)\frake_{\nu'}.
\end{align*}
The sum over $\nu''$ is equal to $|A_m|$ if $\mu\in A^m$ and $0$
otherwise. Hence
\begin{align*}
\frake_{\lambda} \mid U \mid \alpha &= \frac{|A_m|}{|A|}
\sum_{\mu\in A^m} \sum_{\substack{\nu\in A^m}} e\left( \mu^2/2+(\mu,\lambda-\nu/m)\right)\frake_{\nu}\\
&=\frac{1}{|A|} \sum_{\mu\in A} \sum_{\substack{\nu\in A^m}} e\left(
(m\mu)^2/2+(m\mu,\lambda-\nu/m)\right)\frake_{\nu}.
\end{align*}

On the other hand, the right hand side of \eqref{toshow} is equal to
\begin{align*}
\frake_{\lambda} \mid \alpha \mid U^{m^2}
&=\frac{1}{|A|}\sum_{\mu,\nu\in A} e\left(
m^2\mu^2/2+(\mu,m\lambda-\nu)\right)
\frake_{\nu}\\
&=\frac{1}{|A||A_m|}\sum_{\mu,\nu\in A} \sum_{\mu'\in A_m}e\left(
m^2\mu^2/2+(\mu+\mu',m\lambda-\nu)\right) \frake_{\nu}.
\end{align*}
The sum over $\mu'$ is equal to $|A_m|$ if $\nu\in A^m$ and $0$
otherwise. Hence
\begin{align*}
\frake_{\lambda} \mid \alpha \mid U^{m^2}
&=\frac{1}{|A|}\sum_{\mu\in A} \sum_{\nu\in A^m} e\left(
m^2\mu^2/2+(\mu,m\lambda-\nu)\right)
\frake_{\nu}\\
&=\frac{1}{|A|}\sum_{\mu\in A} \sum_{\nu\in A^m} e\left(
(m\mu)^2/2+(m\mu,\lambda-\nu/m)\right) \frake_{\nu}.
\end{align*}
This concludes the proof of the proposition.
\end{proof}

\begin{lemma}\label{lem:cosetact}
Let $\delta \in \tilde\Gamma \alpha \tilde \Gamma$, and let
$\gamma,\gamma'\in\tilde\Gamma$. Then
\[
\frake_\lambda  \mid (\gamma\delta\gamma') = \frake_\lambda  \mid
\gamma \mid \delta\mid \gamma'.
\]
\end{lemma}

\begin{proof}
This follows immediately from the definition and Proposition
\ref{welldef}.
\end{proof}

The element $\beta=\left(
\kzxz{1}{0}{0}{m^2},m \right)\in \widetilde{\Gl}_2^+(\R)$ belongs to the double coset $\tilde \Gamma \alpha \tilde \Gamma$. The following Proposition gives its action on $\C[A]$.

\begin{proposition}\label{prop:adj}
We have
\begin{align}
\frake_\lambda \mid \beta = \sum_{\substack{\mu\in A\\ m\mu = \lambda}}
\frake_{\mu}.
\end{align}
Moreover, for the standard scalar product on $\C[A]$ we have
\begin{align}
\langle \fraka \mid \alpha , \frakb\rangle = \langle \fraka  , \frakb\mid \beta\rangle,\qquad \fraka,\frakb\in\C[A].
\end{align}
\end{proposition}

\begin{proof}
The first assertion follows from the fact that $\beta=S\alpha S^{-1}$ and Lemma
\ref{lem:cosetact}.
The second assertion is verified by a straightforward computation.
\end{proof}

\begin{proposition}\label{prop:product}
Let $m, n$ be coprime positive integers, and put  $\alpha = \left(\kzxz{m^2}{0}{0}{1},
  1\right)$ and  $\beta = \left(\kzxz{n^2}{0}{0}{1},1\right)$. Then for  $g\in \tilde\Gamma\alpha\tilde\Gamma$ and $h\in
\tilde\Gamma\beta\tilde\Gamma$ we have
\begin{align}\label{eq:act}
\frake_\lambda\mid g\mid h = \frake_\lambda\mid (gh).
\end{align}
\end{proposition}

\begin{proof}
We write $g=\gamma \alpha \gamma'$ and $h=\delta\beta\delta'$ with $\gamma,\gamma',\delta,\delta'\in \tilde\Gamma$. Since $(m,n)=1$, a simple argument using the elementary divisor theorem shows that  $gh=\epsilon \alpha \beta \epsilon'$ for suitable $\epsilon,\epsilon'\in \tilde\Gamma$.

In view of Lemma \ref{lem:cosetact} it suffices to prove the assertion in the case that $\gamma=\delta'=1$.

As a second reduction step, we now show that we may in addition assume that $\epsilon'=1$. In fact, the identity $gh=\alpha\gamma'\delta\beta=
\epsilon\alpha\beta\epsilon'$ implies that
\[
\delta\beta\epsilon'{}^{-1}=\gamma'{}^{-1} \alpha^{-1}\epsilon\alpha\beta.
\]
The matrix component of the left hand side has integral entries, hence the same is true for the right had side. Using the coprimality of $m$ and $n$ we may infer that
\[
\tilde\delta:=\gamma'{}^{-1} \alpha^{-1}\epsilon\alpha
\]
belongs to $\tilde\Gamma$.
We obtain
\begin{align*}
\delta\beta&=\tilde\delta\beta\epsilon',\\
(\alpha\gamma')(\tilde\delta\beta)&=\epsilon\alpha\beta.
\end{align*}
From the assertion in the $\epsilon'=1$ case we get
\[
\frake_\lambda\mid (\alpha\gamma')\mid (\tilde\delta\beta)= \frake_\lambda\mid (\epsilon\alpha\beta),
\]
which implies the assertion for arbitrary $\epsilon'$.

Finally we need to prove the claim in the case that
$\gamma=\delta'=\epsilon'=1$. So $g=\alpha \gamma'$, $h=\delta\beta$ and
\[
gh=\alpha \gamma'\delta\beta=\epsilon\alpha\beta.
\]
Since $\alpha \gamma'\delta=\epsilon\alpha$, Lemma \ref{lem:cosetact} implies
that
\begin{align*}
\frake_\lambda\mid(\alpha \gamma'\delta)&=\frake_\lambda\mid(\epsilon\alpha),\\
\frake_\lambda\mid(\alpha \gamma'\delta)\mid\beta&=\frake_\lambda\mid(\epsilon\alpha)\mid \beta.
\end{align*}
Now the claim follows from Lemma \ref{lem:cosetact} and the fact that
$\frake_\lambda \mid \alpha\mid\beta = \frake_\lambda\mid (\alpha\beta )$.
\end{proof}

\begin{definition}
Let $m$ be a positive integer and $\alpha=\left(
\kzxz{m^2}{0}{0}{1},1 \right)\in \widetilde{\Gl}_2^+(\R)$. Let
\[
\tilde \Gamma \cdot \alpha \cdot \tilde \Gamma = \bigcup_{i}
\tilde\Gamma \cdot \delta_i
\]
be a disjoint left coset decomposition. We define the Hecke operator
$T(m^2)^*$ on modular forms  $f\in M_{k,A}$  by
\[
f\mapsto f\mid_{k,A} T(m^2)^* = m^{k-2} \sum_i f\mid_{k,A} \delta_i
= m^{k-2} \sum_i \sum_{\lambda\in A} \big(f_\lambda \mid_{k}
\delta_i\big)\otimes \big(\frake_\lambda \mid_A \delta_i\big).
\]
\end{definition}

Lemma \ref{lem:cosetact} implies that the definition does not depend
on the choice of the coset representatives.  Notice that for $m$
coprime to $N$, this definition agrees with the earlier definition
in Section \ref{sect:4.3}.


\begin{theorem}
\label{thm:structure2}
For any positive integer $m$, the Hecke operator $T(m^2)^*$ is a linear operator on $M_{k,A}$
taking cusp forms to cusp forms. It is self adjoint with respect to
the Petersson scalar product. Moreover, if $m,n$ are coprime, then
\[
T(m^2)^*T(n^2)^* = T(m^2 n^2)^*.
\]
\end{theorem}

\begin{proof}
  The first assertion is a consequence of Lemma \ref{lem:cosetact}.
  The self adjointness follows from Proposition \ref{prop:adj} along
  the standard argument for scalar valued modular forms
  (cf.~\cite{Bu}, Theorem 1.4.3). The last assertion is a consequence
  of Proposition \ref{prop:product} and the corresponding property of
  the abstract Hecke algebra.
\end{proof}

\begin{remark}
For a prime $p$ dividing $N$ the local Hecke algebra, that is, the
subalgebra of $\End(M_{k,A})$ generated by the $T(p^{2\nu})^*$, is
considerably more complicated than in the case where $p$ is coprime
to $N$. For instance, it is commutative if $p$ is coprime to $N$,
but in general non-commutative if $p$ divides $N$.
\end{remark}


\end{document}